\documentclass[12pt,twoside,reqno]{amsart}
\usepackage{amssymb}
\usepackage{amsmath}
\usepackage{amsbsy}
\font\teneufm=eufm10 scaled \magstep1
\font\seveneufm=eufm7 scaled \magstep1
\font\fiveeufm=eufm5  scaled \magstep1
\newfam\eufmfam
\textfont\eufmfam=\teneufm
\scriptfont\eufmfam=\seveneufm
\scriptscriptfont\eufmfam=\fiveeufm
\def\frak#1{{\fam\eufmfam\relax#1}}

\usepackage{mathrsfs}

\newfam\msbfam
\font\tenmsb=msbm10 scaled \magstep1  \textfont\msbfam=\tenmsb
\font\sevenmsb=msbm7 scaled \magstep1 \scriptfont\msbfam=\sevenmsb
\font\fivemsb=msbm5 scaled \magstep1  \scriptscriptfont\msbfam=\fivemsb

\def\One{{1\kern-3.8pt 1}}
\def\one{{1\kern-3.1pt 1}}

\def\RR{{\mathbb R}}
\def\CC{{\mathbb C}}

\def\ZZ{{\mathbb Z}}

\def\PP{{\mathbb P}}

\def\GG{{\mathbb G}}

\def\dd#1{\raise1.5pt\hbox{$\,\partial\!$}/\raise-2.5pt\hbox{$\!\partial#1\,$}}

\def\reg{\mathop{\rm reg}\nolimits}
\def\Bir{\mathop{\rm Bir}\nolimits}
\def\GL{\mathop{\rm GL}\nolimits}
\def\PGL{\mathop{\rm PGL}\nolimits}

\def\Aut{\mathop{\rm Aut}\nolimits}
\def\Im{\mathop{\rm Im}\nolimits}

\def\ad{\mathop{\rm ad}\nolimits}
\def\End{\mathop{\rm End}\nolimits}

\def\BR{\mathop{\rm BR}\nolimits}

\def\ra{\rightarrow}

 \def\HollowBoxx #1#2#3{{\dimen0=#1 \advance\dimen0 by -#2
       \dimen1=#1 \advance\dimen1 by #3
        \vrule height 0pt depth #3 width #2
       \hskip -#3
       \vrule height #1 depth #3 width #3}}
 \def\LeftContraction{\mathord{\kern1.45pt \HollowBoxx{6pt}{3.5pt}{.4pt}}\,}

 \def\HollowBox #1#2#3{{\dimen0=#1 \advance\dimen0 by -#3
       \dimen1=#1 \advance\dimen1 by #3
        \vrule height #1 depth #3 width #3
        \vrule height 0pt depth #3 width #2
        \hskip -#3}}
 \def\RightContraction{\mathord{\, \HollowBox{6pt}{3.1pt}{.4pt}} \kern1.6pt}

\def\qed{{\hfill $\Box$}}
\newtheorem{theorem}{THEOREM}[section]
\newtheorem{corollary}[theorem]{Corollary}

\newtheorem{example}[theorem]{Example}
\newtheorem{remark}[theorem]{Remark}
\newtheorem{proposition}[theorem]{Proposition}

\newtheorem{definition}[theorem]{Definition}

\begin{document}

\begin{center}
{\large \bf Regularization of Local CR-Automorphisms 
\vspace{0.3cm}\\
of Real-Analytic CR-Manifolds}
\vspace{0.5cm}\\
\normalsize A. Isaev and W. Kaup
\end{center}

\thispagestyle{empty}

\pagestyle{myheadings}
\markboth{A. Isaev and W. Kaup}{Local CR-Automorphisms of Real-Analytic CR-Manifolds}

\begin{quotation}\small {\bf Abstract:} Let $M$ be a connected generic real-analytic CR-sub\-man\-i\-fold of a finite-dimensional complex vector space $E$. Suppose that for every $a\in M$ the Lie algebra ${\frak{hol}}(M,a)$ of germs of all infinitesimal real-analytic CR-automorphisms of $M$ at $a$ is finite-dimensional and its complexification contains all constant vector fields $\alpha\dd z$, $\alpha\in E$, and the Euler vector field $z\dd z$. Under these assumptions we show that: (I) every ${\frak{hol}}(M,a)$ consists of polynomial vector fields, hence coincides with the Lie algebra ${\frak{hol}}(M)$ of all infinitesimal real-analytic CR-automorphisms of $M$; (II) every local real-analytic CR-auto\-mor\-phism of $M$ extends to a birational transformation of $E$, and (III) the group $\Bir(M)$ generated by such birational transformations is realized as a group of projective transformations upon embedding $E$ as a Zariski open subset into a projective algebraic variety. Under additional assumptions the group $\Bir(M)$ is shown to have the structure of a Lie group with at most countably many connected components and Lie algebra ${\frak{hol}}(M)$. All of the above results apply, for instance, to Levi non-degenerate quadrics, as well as a large number of Levi degenerate tube manifolds.        
\end{quotation}

\setcounter{section}{0}

\section{Introduction and Preliminaries}\label{intro}
\setcounter{equation}{0}

Let $h=(h_1,\dots,h_k)$ be a $\CC^k$-valued Hermitian form on $\CC^n$, with $n,k\ge 1$. The form $h$ is called {\it non-degenerate}\, if the following two conditions are satisfied:
\vspace{-0.3cm}\\

\noindent (i) the scalar Hermitian forms $h_1,\dots,h_k$ are
linearly independent\linebreak over $\RR$;
\vspace{-0.3cm}\\

\noindent (ii) $h(z,z')=0$ for all $z'\in\CC^n$ implies $z=0$.
\vspace{-0.3cm}\\

\noindent For a non-degenerate $h$ one has $k\le n^2$. Note that many authors define a non-degenerate Hermitian form as a form satisfying condition (ii) alone. 

To any $\CC^k$-valued Hermitian form $h$ on $\CC^n$ one associates the quadric $Q_h\subset\CC^{n+k}$ of CR-dimension $n$ and CR-codimension
$k$ as follows:
$$
Q_h:=\{(z,w)\in\CC^{n+k}:\Im w=h(z,z)\},
$$
where $z=(z_1,\dots,z_n)$ is a point in $\CC^n$, and $w=(w_1,\dots,w_k)$ is a point in $\CC^k$. The CR-manifold $Q_h$ is called the {\it quadric associated to $h$}. 

If $h$ is non-degenerate, then any $C^1$-smooth CR-isomorphism between domains in $Q_h$ extends to a birational map of $\CC^{n+k}$ (see the classical papers \cite{Po}, \cite{Tan1}, \cite{A} for $k=1$ and the papers \cite{KT}, \cite{F}, \cite{Tum}, \cite{Ka1}, \cite{Su}, \cite{B1}, \cite{B2} for $1<k\le n^2$). These birational maps form a group (this is not obvious at all and requires a justification -- see Remark \ref{propstar}). We denote this group by $\Bir(Q_h)$ and call the {\it group of birational transformations of $Q_h$}. For $k=1$ every element of $\Bir(Q_h)$ is a linear fractional transformation induced by an automorphism of $\CC\PP^{n+1}$ (see \cite{Po}, \cite{Tan1}, \cite{A}). For some Hermitian forms $h$ with $1<k\le n^2$ formulas for elements of certain subgroups of $\Bir(Q_h)$ were given in \cite{ES2}, \cite{ES3}. It was shown in \cite{Tum} that the group $\Bir(Q_h)$ can be endowed with the structure of a Lie group (possibly with uncountably many connected components) with Lie algebra isomorphic to the Lie algebra of all infinitesimal CR-automorphisms of $Q_h$, where a smooth vector field on $Q_h$ is called an infinitesimal CR-automorphism if in a neighborhood of every point of $Q_h$ translations along the integral curves of the vector field form a local group of local CR-automorphisms. Every infinitesimal CR-automorphism of $Q_h$ is known to be polynomial. We will see below that $\Bir(Q_h)$ can be embedded in a natural way into the complex group $\PGL_N(\CC)$ as a closed real subgroup (see Corollary \ref{generalcor} and Remark \ref{quadr}). 

We are interested in regularizing the elements of the group $\Bir(Q_h)$ as stated in Definition \ref{regulariz} below. This definition applies to more general CR-submanifolds $M$ of a finite-dimensional complex vector space $E$ than quadrics, and we will first introduce $\Bir(M)$, the {\it group of birational transformation of $M$}. Throughout the paper $M$ is assumed to be connected, locally closed, real-analytic and generic in $E$.

For a rational map $g$ of $E$, we denote by $\reg(g)$ the subset of all regular points of $g$. Let $\Bir(E)$ be the group of all birational transformations of $E$. The restriction of $g\in\Bir(E)$ to the Zariski open subset $\reg(g)$ defines a biholomorphic map $\reg(g)\ra\reg(g^{-1})$. Denote by $\BR(M)$ the collection of all elements $g\in\Bir(E)$ with the following property: there exists a non-empty domain $V\subset M$ with $V\subset\reg(g)$ and $g(V)\subset M$. Note that in general elements of $\BR(M)$ are not defined on all of $M$: they may have poles and points of indeterminacy on $M$. It is clear that $(\BR(M))^{-1}=\BR(M)$. In general, however, $\BR(M)$ is a proper subset of $\BR(M)\cdot\BR(M)$. We define $\Bir(M)$ to be the subgroup of $\Bir(E)$ generated by $\BR(M)$. 

One can give a sufficient condition that guarantees that $\Bir(M)=\BR(M)$. Recall, first of all, that $M$ is called {\it minimal at a point}\, $a\in M$, if there does not exist a CR-submanifold $M_0\subset M$, with $\dim M_0<\dim M$ and\linebreak CR-dim$M_0$=CR-dim$M$, passing through $a$. The manifold $M$ is called {\it minimal}\, if it is minimal at its every point. 

Let $M$ be a connected real-analytic generic CR-submanifold of $E$. For such $M$ we introduce the following
$$
\begin{array}{l}
\hbox{{\bf Condition \boldmath{$(*)$}:}}\\
\vspace{-0.1cm}\\
\hbox{(a) $M$ is minimal,}\\
\vspace{-0.1cm}\\
\hbox{(b) $M_1\subset M$ holds for every connected real-analytic submanifold}\\
\vspace{-0.3cm}\\
\hspace{0.7cm}\hbox{$M_1\subset E$ such that $W\cap M=W\cap M_1\ne\emptyset$ for some domain $W$ in $E$.}
\end{array}
$$
In Proposition \ref{cond8} in Section \ref{birational} we show that if $M$ satisfies Condition $(*)$ then $\Bir(M)$ coincides with $\BR(M)$. This condition is satisfied, for example, if $M$ is minimal and closed in $E$. In particular, Condition $(*)$ is satisfied for any quadric $Q_h$ (note that part (i) of the definition of the non-degeneracy of an Hermitian form $h$ given at the beginning of Section \ref{intro} is equivalent to $Q_h$ being minimal). There are also a large number of examples of non-closed everywhere Levi degenerate CR-submanifolds satisfying Condition $(*)$. An interesting family of such CR-submanifolds is presented in Example \ref{exx2} in Section \ref{propp}.

\begin{remark}\label{propstar}\rm Proposition \ref{cond8} plays a key role in understanding the group $\Bir(Q_h)$, but it appears to have been overlooked in the literature on quadrics so far. Indeed, many authors seem to assume without proof that the set of maps $\BR(Q_h)$ is a group.
\end{remark}     

We will now give an exact definition of what we mean by regularization. For a complex manifold $Y$ we denote by $\Aut(Y)$ the group of all biholomorphic automorphisms of $Y$.

\begin{definition}\label{regulariz}\sl Let $M$ be a connected real-analytic generic CR-submanifold $M$ of a finite-dimensional complex vector space $E$. A subgroup $G\subset\Bir(M)$ is said to be 
\vspace{0.3cm}

\noindent {\rm (i)} {\it regularizable}\, on a complex manifold $Y$, if there exists an open holomorphic embedding $\varphi:E\ra Y$ and a group homomorphism $\tau: G\ra\Aut(Y)$, such that for every $g\in G$ one has $\varphi\circ g=\tau(g)\circ\varphi$ on $\reg(g)$;
\vspace{0.3cm}

\noindent {\rm (ii)} {\it projectively regularizable} if for a suitable integer $N$ there exists an irreducible complex algebraic subvariety $X\subset\CC\PP^N$, a group homomorphism $\tau: G\ra\PGL_{N+1}(\CC)$, and an algebraic isomorphism $\varphi: E\mapsto X_0$, where $X_0$ is a Zariski open subset of $X$, such that $\varphi\circ g=\tau(g)\circ\varphi$ on $\reg(g)$ for every $g\in G$.
\vspace{0.3cm}

\noindent The map $\varphi$ is called a {\it regularization map}.

\end{definition}

\noindent Clearly, if $G$ is projectively regularizable, it is regularizable on the connected Zariski open subset
\begin{equation}
\hat E:=\bigcup_{g\in G}\tau(g)\varphi(E),\label{setZ}
\end{equation}
of the non-singular part $X_{\reg}$ of $X$. The set $\hat E$ is the smallest $\tau(G)$-invariant domain in $X$ that contains $\varphi(E)$. Note also that one can assume that $X$ is not contained in any projective hyperplane in $\CC\PP^N$.  

Regularization results for certain groups of birational transformations can be found in \cite{HZ}, \cite{Z1}. If $Q_h$ is a hyperquadric (i.e. $k=1$), the group $\Bir(Q_h)$ is known to be projectively regularizable with $N=n+1$ due to the classical work \cite{Po}, \cite{Tan1}, \cite{A}. Further, it was shown in \cite{ES1} (see also \cite{B2}, \cite{Mi}) that for $2\le k\le n^2-1$, excluding the situation $k=n=2$, a quadric in general position has only affine automorphisms, in which case $\Bir(Q_h)$ is projectively regularizable with $N=n+k$ for trivial reasons. In fact, we show in Section \ref{generalproof} that $\Bir(Q_h)$ is projectively regularizable for any non-degenerate form $h$. This is a consequence of our main theorem, which applies to much more general CR-manifolds than quadrics. In order to state the theorem we need to introduce some notation and give necessary definitions.

Let $M$ be a real-analytic generic CR-submanifold of a complex manifold $Z$. In what follows all local CR-automorphisms and infinitesimal CR-automorphism of $M$ are assumed to be {\it real-analytic} (note that a $C^1$-smooth CR-isomorphism between Levi non-degenerate real-analytic CR-manifolds is in fact real-analytic -- see Theorem 3.1 in \cite{BJT}). We denote by ${\mathfrak{hol}}(M)$ the real Lie algebra of all real-analytic infinitesimal CR-automorphisms of $M$. A vector field $\xi$ on $M$ lies in ${\mathfrak{hol}}(M)$ if and only if $\xi$ extends to a holomorphic vector field on a neighborhood $U$ of $M$ in $Z$. [We think of holomorphic vector fields on $U$ as holomorphic sections over $U$ of the tangent bundle $TU$. In particular if $Z=E$, a holomorphic vector field $f(z)\dd z$ is just given by a holomorphic map $f: U\ra E$.] 

For $a\in M$ we denote by ${\mathfrak{hol}}(M,a)$ the real Lie algebra of all germs at $a$ of vector fields in ${\mathfrak{hol}}(V)$, with $V$ running over all open neighborhoods of $a$ in $M$. Clearly, ${\mathfrak{hol}}(M,a)$ is a real Lie subalgebra of the complex Lie algebra ${\mathfrak{hol}}(Z,a)$. By Proposition 12.5.1 of \cite{BER} the finite-dimensionality of ${\mathfrak{hol}}(M,a)$ implies that $M$ is {\it holomorphically non-degenerate at}\, $a$, i.e. the Lie algebra ${\mathfrak{hol}}(M,a)$ is totally real in ${\mathfrak{hol}}(Z,a)$ for all $a\in M$. Indeed, if $\xi$ lies in ${\mathfrak{hol}}(M,a)\cap i{\mathfrak{hol}}(M,a)$, then $\psi\cdot\xi\in{\mathfrak{hol}}(M,a)$ for any germ $\psi$ of a holomorphic function near $a$. Thus the formal complexification of ${\mathfrak{hol}}(M,a)$ is isomorphic to ${\mathfrak{hol}}(M,a)+i{\mathfrak{hol}}(M,a)\subset{\mathfrak{hol}}(Z,a)$ if $\hbox{dim}\,{\mathfrak{hol}}(M,a)<\infty$.

Let $M$ be a real-analytic generic CR-submanifold of a finite-dimensional complex vector space $E$. For such $M$ and a point $a\in M$  we introduce the following 
$$
\begin{array}{l}
\hbox{{\bf Property (P) at $a$:}}\\
\vspace{-0.1cm}\\
\hbox{(a) the Lie algebra ${\mathfrak{hol}}(M,a)$ is finite-dimensional,}\\
\vspace{-0.1cm}\\
\hbox{(b) the complex Lie algebra ${\mathfrak{hol}}(M)+i{\mathfrak{hol}}(M)$ contains the complex solv-}\\
\vspace{-0.3cm}\\
\hspace{0.7cm}\hbox{able Lie algebra}
\end{array}
$$
\begin{equation}
{\mathfrak s}:=\left\{(\alpha+cz)\dd z:\alpha\in E,c\in\CC\right\}.\label{s} 
\end{equation}
Further, we say that $M$ has {\bf Property (P)} if it has Property (P) at its every point. In Section \ref{propp} we give sufficient conditions for $M$ to have Property (P) (see Proposition \ref{suffcondpropp}), and discuss several examples. In particular, every non-degenerate quadric $Q_h$ has Property (P).
 
We now state our main result that provides projective regularization of $\Bir(M)$ for a large class of CR-submanifolds.  

\begin{theorem}\label{general}\sl Let $M$ be a connected real-analytic generic CR-sub\-man\-i\-fold of a finite-dimensional complex vector space $E$. Assume further that $M$ has Property {\rm (P)}. Then the following holds:
\vspace{-0.3cm}\\

\noindent{\rm (I)}  for every $a\in M$ the Lie algebra ${\mathfrak{hol}}(M,a)$ consists of polynomial vector fields, hence ${\mathfrak{hol}}(M,a)={\mathfrak{hol}}(M)$; 
\vspace{-0.3cm}\\ 

\noindent{\rm (II)} every real-analytic CR-isomorphism $g$ between non-empty domains in $M$ extends to a map lying  in $\Bir(M)$ of the form $q(z)^{-1}p(z)$, where $p:E\ra E$, $q:E\ra\End(E)$ are polynomial maps, and $\reg(g)=\reg(q^{-1})=\{z\in E:\det q(z)\ne 0\}$;
\vspace{-0.3cm}\\

\noindent{\rm (III)} $\Bir(M)$ is projectively regularizable.
\end{theorem}

Our next theorem provides information on the extension of $\varphi(M)$ into $\CC\PP^N$. Recall that a real-analytic CR-manifold $M$ is called {\it locally homogeneous at a point}\, $a\in M$ if the evaluation map ${\mathfrak{hol}}(M,a)\ra T_aM$, $\xi\mapsto\xi_a$, is surjective, and $M$ is called {\it locally homogeneous}\, if $M$ is locally homogeneous at every point (see \cite{Z2} for equivalent definitions of local homogeneity). In the theorem to follow we assume that $M$ has Property (P) at some point, satisfies part (b) of Condition $(*)$, and is locally homogeneous. Observe that these assumptions imply that $M$ has Property (P) and satisfies Condition $(*)$. Indeed, local homogeneity implies that $M$ has Property (P). Further, by Proposition 4.2 of \cite{Z2} the finite-dimensionality of ${\mathfrak{hol}}(M,a)$ and local homogeneity at $a$ for all points $a\in M$ yield that $M$ is minimal. Hence $M$ satisfies Condition $(*)$. For such a manifold $M$, we denote by $\hat M$ the unique $\Bir(M)$-orbit in $\CC\PP^N$ containing $\varphi(M)$. Clearly, $\hat M$ is a connected immersed generic CR-submanifold of $\hat E$ (see (\ref{setZ})), and we denote by $\Aut(\hat M)$ the group of all real-analytic CR-automorphisms of $\hat M$.

We now state our next result.      

\begin{theorem}\label{general2}\sl Let $M$ be a connected real-analytic generic CR-sub\-man\-i\-fold of a finite-dimensional complex vector space $E$. Assume that $M$ has Property {\rm (P)} at some point, satisfies part (b) of Condition $(*)$, and is locally homogeneous. Then for the regularization map $\varphi$ and homomorphism $\tau$ arising in Theorem \ref{general} the set $\varphi(M)$ is open and dense in $\hat M$, and $\Aut(\hat M)=\tau(\Bir(M))$. Furthermore, if $\overline{M}\setminus M$ does not contain a CR-submanifold of $E$ locally CR-equivalent to $M$, then $\tau(\Bir(M))$ is closed in $\PGL_{N+1}(\CC)$, and the Lie algebra of $\tau(\Bir(M))$ is canonically isomorphic to ${\mathfrak{hol}}(M)$.
\end{theorem}

For $M$ satisfying the assumptions of Theorem \ref{general2} we now introduce a Lie group structure on $\Bir(M)$ by pulling back the Lie group structure from $\tau(\Bir(M))$ by means of $\tau$. In this Lie group topology $\Bir(M)$ has at most countably many connected components. In Section \ref{liegrstr} we give another sufficient condition for the existence of a Lie group structure on $\Bir(M)$ with this property (see Theorem \ref{liegr}). It comes from the natural faithful representation of $\Bir(M)$ on ${\mathfrak{hol}}(M)$. 

Applying Theorems \ref{general}, \ref{general2}, \ref{liegr} to any quadric $Q_h$ we obtain the following corollary.

\begin{corollary}\label{generalcor}\sl If $h$ is non-degenerate, then $\Bir(Q_h)$ is projectively regularizable, and for the regularization map $\varphi$ the set $\varphi(Q_h)$ is open and dense in a $\Bir(Q_h)$-orbit in $\CC\PP^N$. The corresponding homomorphism $\tau$ maps $\Bir(Q_h)$ onto a closed real subgroup of\, $\PGL_{N+1}(\CC)$, and $\Bir(Q_h)$ admits the structure of a Lie group with at most countably many connected components and Lie algebra isomorphic to ${\mathfrak{hol}}(Q_h)$. 
\end{corollary}

For the case when $Q_h$ is the \v Silov boundary of a Siegel domain, the regularization statement of Corollary \ref{generalcor} is essentially contained in Theorem 9 of \cite{KMO}.

\begin{remark}\label{quadr}\rm For quadrics the degrees of the polynomial maps $p$ and $q$ arising in statement (II) of Theorem \ref{general} do not exceed 2. The rationality property for local automorphisms of quadrics can be derived already from the results of \cite{Ka1} (see Satz 2, p. 134). This property was also obtained in \cite{Tum}, but our arguments are simpler even for more general CR-manifolds. In addition, a Lie group structure on $\Bir(Q_h)$ with Lie algebra ${\mathfrak{hol}}(Q_h)$ has been constructed in \cite{Tum} by means of considering the natural faithful representation $\rho$ of $\Bir(Q_h)$ on ${\mathfrak{hol}}(Q_h)$ that maps every $g\in\Bir(Q_h)$ into the corresponding push-forward transformation $g_*$ of vector fields in ${\mathfrak{hol}}(Q_h)$. By a general theorem due to Palais (see \cite{Pa}, Theorem VII, p. 103) the image $\rho(\Bir(Q_h))\subset\GL({\mathfrak{hol}}(Q_h))$ has the structure of a Lie group with Lie algebra ${\mathfrak{hol}}(Q_h)$, but this Lie group may a priori have uncountably many connected components, if $\rho(\Bir(Q_h))$ is not closed in $\GL({\mathfrak{hol}}(Q_h))$. No proof of closedness was given in \cite{Tum}. Our construction of a Lie group structure on $\Bir(M)$ in Theorem \ref{general2} relies on the algebraic regularization map $\varphi:E\ra\CC\PP^N$, while the Lie group structure arising in Theorem \ref{liegr} comes from the natural representation $\rho$ of $\Bir(M)$ on ${\mathfrak{hol}}(M)$. In Theorem \ref{general2} we show that $\Bir(M)$ embeds as a closed subgroup into $\PGL_{N+1}(\CC)$, whereas in Theorem \ref{liegr} we prove that $\rho(\Bir(M))$ is closed in $\GL({\mathfrak{hol}}(M))$. The Lie group structures on $\Bir(Q_h)$ arising from Theorems \ref{general2} and \ref{liegr} for $M=Q_h$ are identical. We also note that since the extension $\hat Q_h$ of $Q_h$ is Levi non-degenerate and has pairwise equivalent Levi forms at all points, the existence of the structure of a Lie group on $\Aut(\hat Q_h)$ (and hence on $\Bir(Q_h)$) with Lie algebra ${\mathfrak{hol}}(Q_h)$ in a certain topology follows from the results of \cite{Tan2}. We refer the reader to \cite{BRWZ}, \cite{LMZ} and references therein for results on the existence of Lie group structures on the groups of CR-automorphisms of more general CR-manifolds. 
\end{remark}

If one does not insist on finding a projective regularization, the group $\Bir(Q_h)$ (in fact, the group $\Bir(M)$ for much more general $M$) can be regularized on some complex manifold in the sense of part (i) of Definition \ref{regulariz}, as follows. Consider the complexification ${\mathfrak l}$ of ${\mathfrak{hol}}(Q_h)$. The complex Lie algebra ${\mathfrak l}$ consists of polynomial vector fields of degree not exceeding 2 and has a natural grading ${\mathfrak l}={\mathfrak l}^{-1}\oplus{\mathfrak l}^0\oplus{\mathfrak l}^1$, where the Lie subalgebra ${\mathfrak l}^{-1}$ consists of all constant vector fields on $E$, and all vector fields in ${\mathfrak l}_0:={\mathfrak l}^0\oplus{\mathfrak l}^1$ vanish at the origin (see e.g. Section \ref{generalproof}). Since $[\xi,{\mathfrak l}^0]$ is not contained in ${\mathfrak l}^0$ for every non-zero $\xi\in{\mathfrak l}^{-1}$, the normalizer of ${\mathfrak l}_0$ in ${\mathfrak l}$ coincides with ${\mathfrak l}_0$. Let ${\mathfrak L}$ be the connected simply-connected group with Lie algebra ${\mathfrak l}$. The stabilizer ${\mathfrak L}_0$ of ${\mathfrak l}_0$ under the adjoint representation of ${\mathfrak L}$ is a closed complex subgroup of ${\mathfrak L}$. Since the normalizer of ${\mathfrak l}_0$ in ${\mathfrak l}$ coincides with ${\mathfrak l}_0$, the Lie algebra of ${\mathfrak L}_0$ coincides with ${\mathfrak l}_0$. Thus ${\mathfrak L}_0^{\circ}$ is a closed complex connected subgroup of ${\mathfrak L}$ with Lie algebra ${\mathfrak l}_0$, and we consider the simply-connected complex homogeneous manifold $Y_h:={\mathfrak L}/{\mathfrak L}_0^{\circ}$. One can show that the vector group $E^{+}:=(E,+)$ naturally lies in ${\mathfrak L}$, and therefore $E$ embeds into $Y_h$ as an an open (and dense) subset. Let $\Bir(Q_h)^{\circ}$ denote the connected component of the identity of $\Bir(Q_h)$ with respect to the Lie group topology on $\Bir(Q_h)$ provided, say, by the results of \cite{Tum}. It can be easily shown that $\Bir(Q_h)^{\circ}$ is regularizable on the manifold $Y_h$.

Further, let $\Bir_0(Q_h):=\{g\in\Bir(Q_h):\, \hbox{$0\in\reg(g)$ and $g(0)=0$}\}$. The full group $\Bir(Q_h)$ is generated by $\Bir(Q_h)^{\circ}$ and $\Bir_0(Q_h)$. For an element $g\in\Bir_0(Q_h)$ the corresponding push-forward map $g_*$ is a Lie algebra automorphism of ${\mathfrak l}$ leaving ${\mathfrak l}_0$ invariant. This automorphism induces an automorphism of ${\mathfrak L}$ leaving ${\mathfrak L}_0^{\circ}$ invariant, and therefore gives rise to an element of $\Aut(Y_h)$. Hence the full group $\Bir(Q_h)$ is regularizable on $Y_h$.

While the approach that we have just outlined solves the regularization problem for $\Bir(Q_h)$ in principle (in the sense of part (i) of Definition \ref{regulariz}), our Theorem \ref{general} contains a much stronger result. It provides an algebraic solution to this problem and applies to a large class of CR-manifolds.   

We would like to thank Michael Eastwood for many helpful discussions. This research is supported by the Australian Research Council and was initiated while the second author was visiting the Australian National University.

\section{Birational Transformations of a Vector Space}\label{birational}
\setcounter{equation}{0}

In this section we state two general propositions on birational maps of a finite-dimensional complex vector space $E$. 

The first proposition will be used in the proofs of Theorems \ref{general}, \ref{general2} but is also of independent interest (cf. \cite{Ko1}, \cite{Ko2}, \cite{Ka2}). For every $\alpha\in E$ we consider the constant holomorphic vector field $\alpha\dd z$, and denote by $\eta$ the {\it Euler vector field}\, $z\dd z$. 
  
\begin{proposition}\label{birat}\sl Let $D_1,D_2\subset E$ be non-empty domains and $g:D_1\to D_2$ a biholomorphic map with induced Lie algebra isomorphism $g_*:{\mathfrak{hol}}(D_1)\ra{\mathfrak{hol}}(D_2)$. With $g^*:=g_*^{-1}$ define the holomorphic maps 
$$
p_g:D_1\ra E\quad \hbox{and}\quad  q_g:D_1\ra\End(E)
$$
by
\begin{equation}
g^*(\eta)=p_g(z)\dd z,\quad g^*(\alpha\dd z)=\big(q_g(z)\alpha\big)\dd z\label{pandq}
\end{equation}
for all $\alpha\in E$. Then $q_g( D_1)\subset\GL(E)$ and \begin{equation}
\hbox{$g(z)=q_g(z)^{-1}p_g(z)$, with $g'(z)=q_g(z)^{-1}$, for all $z\in D_1$.}\label{explform}
\end{equation} 
\end{proposition}

\noindent {\bf Proof:} For every $h(z)\dd z\in{\mathfrak{hol}}(D_2)$ we have by definition
$$
g^*(h(z)\dd z)=\Bigl(g'(z)^{-1}h(g(z))\Bigr)\dd z\in{\mathfrak{hol}}(D_1),
$$
where $g'(z)\in\GL(E)$ for $z\in D_1$, is the derivative of $g$ at $z$. For $h(z)\equiv \alpha$, with $\alpha\in E$, this implies $g'(z)^{-1}\alpha=q_g(z)\alpha$, and for $h(z)\equiv z$ we get $g'(z)^{-1}g(z)=p_g(z)$. Formula (\ref{explform}) follows from these two relations. \qed
\vspace{0.3cm}\\

Recall that ${\mathfrak s}$ is the complex solvable Lie subalgebra of ${\mathfrak{hol}}(E)$ spanned by all constant vector fields $\alpha\dd z$ and the Euler vector field $\eta$ (see (\ref{s})). Proposition \ref{birat} yields the following corollary.

\begin{proposition}\label{corprop1}\sl Suppose that for the biholomorphic map $g:D_1\ra D_2$ from Proposition \ref{birat} all vector fields in both $g^*({\mathfrak s})$ and $g_*({\mathfrak s})$ extend to rational vector fields on $E$. Then $g$ extends to an element of $\Bir(E)$, with $\reg(g)=\reg(g')=\reg(q^{-1})$.
\end{proposition}

\noindent{\bf Proof:} We only need to show that $\reg(g)=\reg(g')$. Clearly, we have $\reg(g)\subset\reg(g')$. To obtain the opposite inclusion, we suppose that $\reg(g')\setminus\reg(g)$ is non-empty. We let $n:=\dim E$, identify $E$ with $\CC^n$, and write $g$ as $g=(g_1,\dots,g_n)$. Then there exists $j$ such that $A:=\reg(g')\setminus\reg(g_j)$ is non-empty. It then follows that one can find a point $a\in A$, which is not an indeterminacy point of $g_j$, that is, $g_j=r_j/s_j$, where $r_j$ and $s_j$ are polynomials with $r_j(a)\ne 0$, $s_j(a)=0$. Hence for some $k$ the order of vanishing of $s_j \partial r_j/\partial z_k- r_j \partial s_j/\partial z_k$ at $a$ is finite and strictly less than that of $s_j^2$. Therefore, $a$ is not a regular point of $\partial g_j/\partial z_k$, which contradicts our choice of $a$. \qed
\vspace{0.3cm}\\

\begin{remark}\label{regset}\rm We will use Proposition \ref{corprop1} in Section \ref{generalproof} in the case when all vector fields in $g^*({\mathfrak s})$ and $g_*({\mathfrak s})$ extend to polynomial vector fields on $E$. In this situation $\reg(g)=\reg(q_g^{-1})=\{z\in E:\det q_g(z)\ne 0\}$. In fact, $\det q_g$ is a denominator of the rational map $g$, that is, $(\det q_g) g$ is a polynomial map. As the following example shows, $\det q_g$ need not be an exact denominator (a denominator of minimal degree) of $g$.
\end{remark} 

\begin{example}\label{exrational}\rm Let $E:=\CC^{n\times m}$, $b\in\CC^{m\times n}$ a fixed
matrix, and 
$$
g(z):=(\One-zb)^{-1}z,
$$
where $\One$ is the $n\times n$ identity matrix. Then $g\in\Bir(E)$ (indeed $g^{-1}(w)=(\One+wb)^{-1}w$). Differentiation yields
$$
g'(z)\alpha=(\One-zb)^{-1}\alpha(\One-bz)^{-1}
$$ 
for all $\alpha\in E$. In particular, for the functions $p,q$ from Proposition \ref{birat} we have
$$
\hbox{$q_g(z)\alpha=(\One-zb)\alpha(\One-bz)$ and $p_g(z)=z-zbz$}
$$
for all $\alpha\in E$. Thus $\det q_g$ is not an exact denominator of $g$. Further, a moment's thought gives $\det q_g(z)=\det(\One-zb)^m\det(\One-bz)^n$, hence $\reg(g)=\{z\in E:\det(\One-zb)\ne 0\}$.
\end{example}
\vspace{0.3cm}

In the next proposition we relate the group $\Bir(M)$ of birational transformations of a CR-submanifold $M\subset E$ to the subset $\BR(M)\subset\Bir(E)$ by means of Condition $(*)$, as stated in Section \ref{intro}. 

\begin{proposition}\label{cond8}\sl Let $M$ be a connected real-analytic generic CR-sub\-man\-i\-fold of $E$. If Condition $(*)$ is satisfied for $M$, then $\Bir(M)=\BR(M)$. Moreover, for every $g\in\BR(M)$ we have $g(M\cap\reg(g))=M\cap\reg(g^{-1})$.
\end{proposition}

\noindent{\bf Proof:} Fix $g\in\BR(M)$, and let $V\subset M$ be a non-empty domain such that $V\subset\reg(g)$ and $g(V)\subset M$. By Lemma 2.2 of \cite{FK2} the non-empty set $M\cap\reg(g)$ is connected, and therefore $M_1:=g(M\cap\reg(g))$ is a real-analytic connected submanifold of $E$. Since $W:=g(V)$ is a non-empty domain in $M$ such that $W\cap M_1=W$, Condition $(*)$ implies that $M_1\subset M\cap\reg(g^{-1})$. Interchanging the roles of $g$ and $g^{-1}$ gives $g(M\cap\reg(g))=M\cap\reg(g^{-1})$.

Now for any $g_1,g_2\in\BR(M)$ we choose a non-empty domain $V\subset M$ with $V\subset\reg(g_1)$ and $g_1(V)\subset\reg(g_2)$. Then $g_2\circ g_1\in\BR(M)$. Therefore, $\BR(M)=\Bir(M)$, as required.\qed
\vspace{0.3cm}

We stress the importance of Proposition \ref{cond8} for the correct understanding of $\BR(M)$ and $\Bir(M)$. In particular, if $M$ does not satisfy the assumptions of Proposition \ref{cond8}, then the set $\BR(M)$ may not be a group. 

As we stated in Section \ref{intro}, a connected real-analytic generic submanifold $M\subset E$ satisfies Condition $(*)$ if $M$ is minimal and closed. There is, however, a large class of examples of non-closed CR-submanifolds satisfying Condition $(*)$. An interesting family of such manifolds is given in Example \ref{exx2} in Section \ref{propp}.

\section{Proof of Theorems \ref{general} and \ref{general2}}\label{generalproof}
\setcounter{equation}{0}

We will first prove Theorem \ref{general}.

Without loss of generality we assume that $M$ contains the origin, and let ${\mathfrak l}$ be the complexification of ${\mathfrak{hol}}(M,0)$. Arguing as in the proof of Proposition 4.2 of \cite{FK1}, we obtain that ${\mathfrak l}$ admits a $\ZZ$-grading 
\begin{equation}
{\mathfrak l}=\bigoplus_{m\in\ZZ}{\mathfrak l}^m,\quad[{\mathfrak l}^m,{\mathfrak l}^{\ell}]\subset{\mathfrak l}^{m+\ell},\label{grad}
\end{equation}
where ${\mathfrak l}^m$ is the $m$-eigenspace of $\ad\eta $ in ${\mathfrak l}$, and ${\mathfrak l}^m=0$ for $m<-1$, as well as for $m$ large enough. Every ${\mathfrak l}^m$ consists of polynomial vector fields homogeneous of degree $m+1$, with 
$$
{\mathfrak l}^{-1}=\{\alpha\dd z:\alpha\in E\}
$$ 
being the Lie algebra of all constant vector fields on $E$. Thus every vector field in ${\mathfrak{hol}}(M,0)$ is polynomial. Arguing in this way for every $a\in M$ we see that all Lie algebras ${\mathfrak{hol}}(M,a)$ are polynomial and hence coincide with ${\mathfrak{hol}}(M)$. Thus we have obtained statement (I). 

For a non-empty domain $D\subset E$ we identify ${\mathfrak l}$ with a Lie subalgebra of ${\mathfrak{hol}}(D)$ by restriction. Let $V_1$, $V_2$ be non-empty domains in $M$, and $g:V_1\to V_2$ a real-analytic CR-isomorphism. Then there exist domains $D_1,D_2\subset E$ and a biholomorphic extension $g:D_1\to D_2$ with $g_*({\mathfrak l})={\mathfrak l}$. Since all vector fields in ${\mathfrak l}$ are polynomial, Proposition \ref{corprop1} yields that $g$ extends to an element of $\Bir(M)$ of the form $q^{-1}p$, where $p:E\ra E$ and $q:E\ra\End(E)$ are polynomial maps (see (\ref{pandq})). By Remark \ref{regset} we have $\reg(g)=\reg(q^{-1})=\{z\in E:\det q(z)\ne 0\}$. Thus we have obtained statement (II).

Further, for every $a\in E$ the isotropy Lie subalgebra 
$$
{\mathfrak l}_a:=\{\xi\in{\mathfrak l}:\xi_a=0\}
$$
has codimension $n:=\dim E$ in ${\mathfrak l}$, and ${\mathfrak l}$ is the direct sum of subspaces ${\mathfrak l}={\mathfrak l}^{-1}\oplus\,{\mathfrak l}_a$, with ${\mathfrak l}_a\ne{\mathfrak l}_b$ for all $a,b\in E$, $a\ne b$. Let $\GG$ be the Grassmannian of all complex linear subspaces $\Lambda\subset{\mathfrak l}$ of codimension $n$. Then $\GG$ is a rational projective algebraic complex manifold on which the complex linear group $\GL({\mathfrak l})$ acts transitively and algebraically by means of the canonical projection $\GL({\mathfrak l})\ra\PGL({\mathfrak l})\subset\Aut(\GG)$. 

The subset 
$$
U:=\{\Lambda\in\GG:{\mathfrak l}={\mathfrak l}^{-1}\oplus\,\Lambda\}
$$ 
is Zariski open in $\GG$ and is algebraically equivalent to the complex vector space of all linear operators $\lambda:{\mathfrak l}_0\to{\mathfrak l}^{-1}$ (just identify every $\lambda$ with its graph $\{\xi+\lambda(\xi):\xi\in{\mathfrak l}_0\}\in\GG$). In this coordinate chart every automorphism of $\GG$ arising from the action of $\GL({\mathfrak l})$ can be written as a matrix linear fractional transformation.

Consider the injective holomorphic map
$$
\varphi:E\to\GG\,,\quad a\mapsto{\mathfrak l}_a.
$$ 
Then $\varphi(E)\subset U$, and since all vector fields in ${\mathfrak l}$ are polynomial, the map $\varphi$ is an algebraic morphism. As a consequence, the set $\varphi(E)$ is constructible. Let $X$ be the Zariski closure of $\varphi(E)$. Clearly, $X$ is an irreducible algebraic subvariety in $\GG$ and $\varphi(E)$ contains a Zariski open (and dense) subset of $X$, hence the topological closure of $\varphi(E)$ in $\GG$ coincides with $X$.

Define
$$
\Bir(E,{\mathfrak l}):=\left\{g\in\Bir(E):g_*({\mathfrak l})={\mathfrak l}\right\}.
$$
Observe that $\Bir(E,{\mathfrak l})$ contains the set $\BR(M)$. Since every element of $\Bir(M)$ is the composition of a finite number of elements of $\BR(M)$, it follows that $\Bir(M)\subset\Bir(E,{\mathfrak l})$.

For any $g\in\Bir(E,{\mathfrak l})$, we regard the push-forward map $g_*$ as an element of $\Aut({\mathfrak l})\subset\GL({\mathfrak l})$, where $\Aut({\mathfrak l})$ is the complex algebraic subgroup of $\GL({\mathfrak l})$ that consists of all Lie algebra automorphisms of ${\mathfrak l}$. Define $\nu$ to be the homomorphism
\begin{equation}
\nu:\,\,\Bir(E,{\mathfrak l})\ra\Aut({\mathfrak l}),\quad g\mapsto g_*.\label{mapPhi}
\end{equation}
By formula (\ref{explform}) the homomorphism $\nu$ is injective. Note that the canonical homomorphism $\pi:\Aut({\mathfrak l})\ra\PGL({\mathfrak l})$ is injective as well. 

Since for $g\in\Bir(E,{\mathfrak l})$ we have $g_*({\mathfrak l}_a)={\mathfrak l}_{g(a)}$ for all $a\in\reg(g)$, the map $\pi(g_*)$ preserves $X$ and the following holds:
\begin{equation}
\varphi\circ g=\sigma(g)\circ\varphi\quad\hbox{on $\reg(g)$,}\label{commutative}
\end{equation}       
where $\sigma:=\pi\circ\nu$. Formula (\ref{commutative}) applies, in particular, to every translation $g(z)=z+\beta$, $\beta\in E$ (note that every translation is an element of $\Bir(E,{\mathfrak l})$). It is straightforward to see that the action of the complex vector group $E^{+}:=(E,+)$ on $\GG$ through the homomorphism $\sigma$ is algebraic, and formula (\ref{commutative}) implies that $\varphi(E)$ is an orbit of this action. It then follows that $\varphi(E)$ is a Zariski open subset of $X$ lying in the non-singular part $X_{\reg}$ of $X$. By Zariski's Main Theorem (see \cite{Mu}, III.9.1), $\varphi:E\to\varphi(E)$ is an algebraic isomorphism. 

We now embed $\GG$ into $\CC\PP^N$ for a sufficiently large integer $N$ by the Pl\"ucker map, and regard $X$ as an algebraic subvariety of $\CC\PP^N$, and $\PGL({\mathfrak l})$ as a subgroup of $\PGL_{N+1}(\CC)$. Formula (\ref{commutative}) then yields that $\Bir(M)$ is projectively regularizable with $\tau:=\sigma|_{\Bir(M)}$. This proves statement (III). 

The proof of Theorem \ref{general} is complete.\qed
\vspace{0.5cm}

We will now prove Theorem \ref{general2}.

Recall that $\hat M$ is a connected generic immersed submanifold of the Zariski open subset $\hat E$ of the non-singular part $X_{\reg}$ of $X$ (see (\ref{setZ})). The minimality of $M$ implies that $\hat M$ is minimal as well. Therefore, it follows from Lemma 2.2 of \cite{FK2} that $\hat M\cap\varphi(E)$ is connected. Since $\hat M$ contains $\varphi(M)$ as an open subset, part (b) of Condition $(*)$ yields that $\hat M\cap\varphi(E)=\varphi(M)$ and that $\varphi(M)$ is dense in $\hat M$.

To show that $\Aut(\hat M)=\tau(\Bir(M))$, observe that for every $g\in \Aut(\hat M)$ there exist domains $V_1, V_2\subset\varphi(M)$ such that $g(V_1)=V_2$. Then the composition $\varphi^{-1}\circ g\circ\varphi$ is a real-analytic CR-diffeomorphism between the domains $\varphi^{-1}(V_1),\varphi^{-1}(V_2)$ in $M$. By statement (II) of Theorem \ref{general}, the map $\varphi^{-1}\circ g\circ\varphi$ extends to an element $g_0$ of $\Bir(M)$, hence $g=\tau(g_0)$. Thus $\Aut(\hat M)=\tau(\Bir(M))$.

Assume now that $\overline{M}\setminus M$ does not contain a CR-submanifold of $E$ locally CR-equivalent to $M$. Let $g_n$ be a sequence in $\tau(\Bir(M))$ converging to an element $g\in\PGL_{N+1}(\CC)$. We claim that $g(\varphi(M))\cap\varphi(E)\ne\emptyset$. Indeed, otherwise $g(\varphi(M))$ lies in a Zariski closed subset of $X_{\reg}$, which is impossible since $g(\varphi(M))$ is generic in $X_{\reg}$. Thus for some domain $V\subset\varphi(M)$ we have $g(V)\subset \varphi(E)$. Clearly, $g(V)$ is a CR-submanifold of $\varphi(E)$ locally equivalent to $\varphi(M)$ and contained in $\overline{\varphi(M)}$. It then follows that there exists $x_0\in V$ for which $g(x_0)\in \varphi(M)$. Since $M$ is locally closed in $E$, for some neighborhood $V'\subset V$ of $x_0$ in $\varphi(M)$ we have $g(V')\subset\varphi(M)$. Hence $g\in\tau(\Bir(M))$, and therefore $\tau(\Bir(M))$ is closed in $\PGL_{N+1}(\CC)$.

Let ${\mathfrak a}$ be the Lie subalgebra of the Lie algebra of $\PGL_{N+1}(\CC)$ corresponding to the closed subgroup $\tau(\Bir(M))\subset \PGL_{N+1}(\CC)$. Every element $v\in{\mathfrak a}$ is a holomorphic vector field on $\CC\PP^N$ tangent to $\hat M$ and gives rise to a holomorphic vector field on $E$ tangent to $M$. 

Conversely, consider a vector field $\xi\in{\mathfrak{hol}}(M)$ and fix $a\in M$. Near $a$ the vector field $\xi$ can be integrated to a local 1-parameter group $t\mapsto g_t$, with $|t|<\varepsilon$, of local real-analytic CR-isomorphisms of $M$. By statement (II) of Theorem \ref{general} every $g_t$ extends to an element of $\Bir(M)$. Further, one can define by composition a map $g_t\in\Bir(M)$ for every $t\in\RR$ and obtain a 1-parameter subgroup of $\Bir(M)$. Then $t\mapsto \tau(g_t)$ is a continuous 1-parameter subgroup of $\tau(\Bir(M))$ and hence $g_t=\exp(tv)$ for some $v\in{\mathfrak a}$.  

Thus we have established an isomorphism between ${\mathfrak a}$ and ${\mathfrak{hol}}(M)$, and the proof of Theorem \ref{general2} is complete.\qed
\vspace{0.3cm}

For the remainder of the article we set ${\mathfrak g}:={\mathfrak{hol}}(M)$. Let $M$ satisfy the assumptions of Theorem \ref{general}, and set $\rho:=\nu|_{\Bir(M)}$, with $\nu$ defined in (\ref{mapPhi}). The homomorphism $\rho$ is injective, and the image $\rho(\Bir(M))$ lies in the group $\Aut({\mathfrak g})$ of all Lie algebra automorphisms of ${\mathfrak g}$, which is a real algebraic subgroup of each of $\GL({\mathfrak g})$ and $\Aut({\mathfrak l})$. Thus $\Bir(M)$ can be endowed with the structure of a Lie group (possibly with uncountably many connected components). An argument identical to that at the end of the proof of Theorem \ref{general2} shows that the Lie algebra of $\Bir(M)$ with respect to this Lie group structure is isomorphic to ${\mathfrak g}$. For $\Bir(M)$ with this Lie group structure the map         
$$
\rho:\Bir(M)\ra\Aut({\mathfrak g})
$$
is just the adjoint representation. In the next section we will show that under additional assumptions  the image $\rho(\Bir(M))$ is closed in $\Aut({\mathfrak g})$.

\section{Closed Embedding of $\Bir(M)$ into $\Aut({\mathfrak g})$}\label{liegrstr}
\setcounter{equation}{0}

Throughout this section we suppose that $M$ has Property (P) and satisfies Condition $(*)$. In particular, by Proposition \ref{cond8} we have $\Bir(M)=\BR(M)$. Under these assumptions, which are weaker than those of Theorem \ref{general2}, we obtain the existence of a Lie group structure on $\Bir(M)$ with at most countably many connected components and Lie algebra ${\mathfrak g}$. Instead of investigating the closedness of $\tau(\Bir(M))$ in $\PGL_{N+1}(\CC)$, as we did in the proof of Theorem \ref{general2}, we investigate the closedness of $\rho(\Bir(M))$ in $\Aut({\mathfrak g})$. Note that a simple example shows that there is always a closed subgroup of $\GL({\mathfrak l})$ whose image in $\PGL({\mathfrak l})$ under the canonical projection $\GL({\mathfrak l})\ra \PGL({\mathfrak l})$ is not closed. 

The result of this section is the following theorem.

\begin{theorem}\label{liegr}\sl Let $M$ be a connected real-analytic generic CR-sub\-man\-i\-fold of a finite-dimensional complex vector space $E$. Assume that $M$ has Property {\rm (P)} and satisfies Condition $(*)$. Then $\rho(\Bir(M))$ is closed in $\Aut({\mathfrak g})$.
\end{theorem}

\noindent {\bf Proof:} Let $g_n$ be a sequence in $\Bir(M)$ such that the sequence $f_n:=\rho(g_n)$ converges in $\Aut({\mathfrak g})$ to an element $f$. By Proposition \ref{birat} every map $g_n$ can be written as $g_n=q_n^{-1}p_n$, where $p_n:=p_{g_n}$ and $q_n:=q_{g_n}$ are polynomial maps on $E$ found from formulas (\ref{pandq}). Since the maps $f_n^{-1}$ converge in $\Aut({\mathfrak g})$ to $f^{-1}$, the sequences $p_n$, $q_n$ converge (uniformly on compact subsets of $E$) to polynomial maps $p:E\ra E$ and $q:E\ra\End(E)$, respectively, and we have $f^{-1}(\alpha\dd z)=\big(q(z)\alpha\big)\dd z$ for all $\alpha\in E$.

Similarly, every map $g_n^{-1}$ can be written as $g_n^{-1}=\tilde q_n^{-1} \tilde p_n$, where $\tilde p_n:=p_{g_n^{-1}}$ and $\tilde q_n:= q_{g_n^{-1}}$ are polynomial maps. The sequences $\tilde p_n$ and $\tilde q_n$ converge to polynomial maps $\tilde p:E\ra E$ and $\tilde q:E\ra\End(E)$, respectively, and we have $f(\alpha\dd z)=\big(\tilde q(z)\alpha\big)\dd z$ for all $\alpha\in E$.

Since $f_k^{-1} f\ra\hbox{id}\in\Aut({\mathfrak g})$, for every neighborhood ${\mathcal V}$ of the identity in $\Aut({\mathfrak g})$ one can find an element $\hat f\in\rho(\Bir(M))$ with $\hat f f\in{\mathcal V}$. Choosing ${\mathcal V}$ such that ${\mathcal V}={\mathcal V}^{-1}$ one can also assume that $f^{-1} \hat f^{-1}\in{\mathcal V}$. Hence by replacing $g_n$ by $\rho^{-1}(\hat f)g_n$ and $f$ by $\hat f f$ we can assume without loss of generality that $\det q\not\equiv 0$ and $\det\tilde q\not\equiv 0$. We then define the rational maps $g:= q^{-1}p$ and $\tilde g:=\tilde q^{-1}\tilde p$.

Let
$$
A:=\{z\in E: \det q(z)=0\},\quad B:=\{z\in E:\det\tilde q(z)=0\}.
$$
Since $\det q_n\ra \det q$ and $\reg(g_n)=\{z\in E:\det q_n(z)\ne 0\}$ (see Proposition \ref{corprop1}), it follows that $g_n\ra g$ uniformly on compact subsets of $E\setminus A$. Similarly, $g_n^{-1}\ra\tilde g$ uniformly on compact subsets of $E\setminus B$. For every $\xi\in{\mathfrak l}$,\linebreak $\xi=h(z)\dd z$, we have
$$
f(\xi)(z)=\tilde q(z)h(\tilde g(z))\dd z,\quad f^{-1}(\xi)(z)=q(z)h(g(z))\dd z.
$$  
Applying the identity
$$
\begin{array}{l}
\xi=f(f^{-1}(\xi))=\tilde q(z)q(\tilde g(z))h(g(\tilde g(z)))\dd z
\end{array}
$$
to constant vector fields $\xi=\alpha\dd z$, we see that $\tilde q(z) q(\tilde g(z))\equiv \hbox{id}$. Applying this identity to the Euler vector field $\eta$ and interchanging the roles of $f$ and $f^{-1}$ we obtain $\tilde g=g^{-1}$, thus $g\in\Bir(E)$. By Proposition \ref{cond8} we have $g_n\in\BR(M)$ and $g_n(M\cap\reg(g_n))=M\cap\reg(g_n^{-1})$, which yields that $g\in\Bir(M)$. Finally, since $g_n\ra g$ on $E\setminus A$ and $g_n^{-1}\ra g^{-1}$ on $E\setminus B$, it follows that $(g_n)_*\ra g_*$. Hence $f=\rho(g)$.

The proof of the theorem is complete. \qed
 
\section{Property (P)}\label{propp}
\setcounter{equation}{0}

In this section we give sufficient conditions for a CR-submanifold to have Property (P) and discuss examples of manifolds with this property (for the statement of Property (P) see Section \ref{intro}).  

We call a CR-submanifold $M$ of a complex manifold $Z$ {\it semi-homogeneous at a point}\, $a\in M$, if the span over $\CC$ of the values $\xi_a$ of elements $\xi\in{\mathfrak{hol}}(M,a)$ at $a$ contains $T_a(M)$. Also, we call $M$ {\it semi-homogeneous}\, if $M$ is semi-homogeneous at every point. Clearly, local homogeneity implies semi-homogeneity for a real-analytic CR-submanifold.

Recall that ${\mathfrak g}={\mathfrak{hol}}(M)$. With this notation we state the following proposition. 

\begin{proposition}\label{suffcondpropp}\sl Let $M$ be a connected real-analytic generic CR-sub\-man\-i\-fold of a finite-dimensional complex vector space $E$. Let $a$ be a point in $M$. Assume that:
$$
\begin{array}{ll}
(1) & \hbox{$M$ is holomorphically non-degenerate at $a$;}\\
\vspace{-0.4cm}\\

(2) & \hbox{$M$ is minimal at $a$;}\\
\vspace{-0.4cm}\\

(3) & \hbox{$M$ is semi-homogeneous at $a$;}\\
\vspace{-0.4cm}\\

(4) & \hbox{the complex Lie algebra ${\mathfrak g}+i{\mathfrak g}$ contains the vector}\\
&\hbox{field $(z-a)\dd z$.}
\end{array}
$$

\noindent Then $M$ has property {\rm (P)}. 
\end{proposition}

\noindent{\bf Proof:} By Theorem 12.5.3 in \cite{BER}, assumptions (1) and (2) imply that ${\mathfrak{hol}}(M,b)$ is finite-dimensional for all $b\in M$.

Without loss of generality we assume that $a=0$. Then using assumption (4) we obtain, as at the beginning of Section \ref{generalproof}, that ${\mathfrak l}:={\mathfrak{hol}}(M,0)+i{\mathfrak{hol}}(M,0)$ admits grading (\ref{grad}) where ${\mathfrak l}^m$ is the $m$-eigenspace of $\ad \eta$ in ${\mathfrak l}$, and ${\mathfrak l}^m=0$ for $m<-1$, as well as for $m$ large enough. Every ${\mathfrak l}^m$ consists of polynomial vector fields homogeneous of degree $m+1$. Since all vector fields in ${\mathfrak l}^m$ for $m>-1$ vanish at the origin, assumption (3) implies that ${\mathfrak l}^{-1}$ is the space of all constant vector fields on $E$. The proof is complete.\qed
\vspace{0.5cm} 

We now give examples of CR-submanifolds that have Property (P).

\begin{example}\label{exx3}\rm Every quadric $Q_h\subset\CC^{n+k}$ associated to a non-degenerate Hermitian form $h$ has Property (P). Indeed, $Q_h$ is homogeneous under the action of the group of maps 
$$
\begin{array}{lll}
z&\mapsto& z+\alpha,\\
\vspace{-0.3cm}\\
w&\mapsto&w+2ih(z,\alpha)+\beta,\\
\end{array}
$$
where $(\alpha,\beta)\in Q_h$. Thus all Lie algebras ${\mathfrak{hol}}(Q_h,a)$ coincide. In fact, they coincide with the Lie algebra ${\mathfrak{hol}}(Q_h)$, which is finite-dimensional (see \cite{B1}, \cite{B2}, \cite{Tum}). Furthermore, ${\mathfrak{hol}}(Q_h)$ clearly contains the vector fields $s\dd w$, $r\dd z+2ih(z,r)\dd z$, $z\dd z+2w\dd w$, $iz\dd z$, with $s\in\RR^k$, $r\in\CC^n$. Therefore, the complexification of ${\mathfrak{hol}}(Q_h)$ contains all constant vector fields and the Euler vector field $\eta$. Hence the quadric $Q_h$ has Property (P). 
\end{example}

\begin{example}\label{exx1}\rm Let $F\subset\RR^{n}$ be an arbitrary connected real-analytic submanifold and
$$
M:=F+i\RR^n\subset \CC^n
$$
the corresponding {\it tube submanifold with base}\, $F$. Then $M$ is a generic semi-homogeneous CR-submanifold of $E$, and ${\mathfrak g}+i{\mathfrak g}$ contains all constant holomorphic vector fields on $E$. Furthermore (see Lemma 4.1 of \cite{FK2}), the tube manifold $M$ is minimal at a point if and only if 
\begin{equation}
\hbox{\sl $F$ is not contained in any affine hyperplane of $\RR^n$,}\label{conda}
\end{equation}
hence a tube manifold is minimal if it is minimal at one point. Next (see Proposition 4.3 of \cite{FK2}), $M$ is holomorphically non-degenerate at a point if and only if
\begin{equation}
\hbox{\sl the only constant vector field $\xi=\alpha\dd x$ tangent to $F$ is $\xi=0$,}\label{condb}
\end{equation}
hence a tube manifold is holomorphically non-degenerate if it is holomorphically non-degenerate at one point. 

To meet conditions (1), (2), (4) of Proposition \ref{suffcondpropp} it is therefore sufficient to require besides (\ref{conda}), (\ref{condb}) that $F$ is a cone, that is, $tF=F$ for every real $t>0$. For every cone $F$ the Levi form of $M$ is degenerate at every point (condition (ii) stated at the beginning of Section \ref{intro} does not hold).
\end{example}

From the large class of tube manifolds that have Property (P) we single out the following special one.

\begin{example}\label{exx2}\rm Fix integers $p\ge q\ge1$ with $n=p+q\ge3$. Then
\begin{equation}
H_{p,q}:=\{x\in\RR^{n}:x^{2}_{1}+\dots
x^{2}_{p}=x^{2}_{p+1}+\dots+x^{2}_{n}\}\label{BA}
\end{equation} 
is a real hyperquadric with 0 as the only singularity. Let $F$ be a connected component of the non-singular part of $H_{p,q}$ and $M:=F+i\RR^{n}$ the corresponding tube submanifold. Then $M$ is a homogeneous CR-submanifold of $\CC^{n}$ that has Property (P) and satisfies Condition $(*)$. The Levi form of $M$ is degenerate at every point. For $q=1$ the non-singular part of $H_{p,q}$ has two connected components (given by $x_n>0$ and $x_n<0$), the future light cone and the past light cone. In this case the group $\Bir(M)$ can be canonically identified with an open subgroup (having two connected components) of $O(n,2)$. For $q>1$ the non-singular part of $H_{p,q}$ is connected, and $\Bir(M)$ is a real algebraic group. For every $a\in M$ the Lie algebra ${\mathfrak{hol}}(M,a)$ is isomorphic to ${\mathfrak{so}}(p{+}1,q{+}1)$ (cf. \cite{FK1}, p. 21).
\end{example}

\begin{example}\label{exx4}\rm Let $D\subset E$ be an irreducible bounded symmetric domain of rank $r$ given in its Harish-Chandra realization, and $Z$ its compact dual containing $E$ as a Zariski open subset. Then $D$ is convex and invariant under the circle group $\exp(i\RR\eta)\subset\GL(E)$. The complex manifold $Z$ is a compact rational algebraic variety on which the simple complex Lie group $L:=\Aut(Z)$ acts transitively. All transformations in $G:=\Aut(D)$ extend to elements of $L$ and, in this way, $G$ is a real form of $L$ and also acts on $Z$. On $Z$ the group $G$ has exactly ${r+2}\choose{2}$ orbits. Among these there are exactly $r+1$ open orbits (including $D$) and a unique closed orbit, the \v Silov boundary of $D$, which coincides with the extremal boundary $\partial_{e}D$ of the convex set $\overline{D}$. Every $G$-orbit is a generic homogeneous CR-submanifold of $Z$ invariant under the circle group $\exp(i\RR\eta)\subset G$. Furthermore, for every orbit $G(b)$, $b\in Z$, the intersection $M:=G(b)\cap E$ is a connected CR-submanifold that has Property (P), and for every $a\in M$ the Lie algebra ${\mathfrak{hol}}(M,a)$ is isomorphic to the Lie algebra of $G$, provided neither $G(b)$ is open in $Z$ nor $G(b)=\partial_{e}D$ in the case when $D$ is of tube type (in this last case $\partial_{e}D$ is totally real). The \v Silov boundary (except when $D$ is of tube type) can be locally realized as a standard quadric in $E$ and, in particular, has non-degenerate Levi form. Every $G$-orbit that is neither open nor closed in $Z$ is Levi degenerate (more precisely 2-nondegenerate) and hence cannot be locally realized as a standard quadric in $E$. The group $\Bir(M)$ is regularizable on the simply-connected complex manifold $Z$, and this is the only possibility up to isomorphism.
\end{example}

\begin{example}\label{exx5}\rm We specialize Example \ref{exx4} to the case
$$
E=\{z\in\CC^{2\times2}:z'=z\}\quad\hbox{and}\quad D=\{z\in E:zz^{*}<\One\},
$$
where $z'$ is the transpose and $z^{*}$ is the transpose conjugate of a matrix $z$. Then $D$ is irreducible symmetric of rank 2, and $Z$ can be identified with a complex projective quadric in $\CC\PP^4$. Also, $G$ is isomorphic to an open subgroup of $O(2,3)$ of index 2. The boundary $\partial D$ of $D$ decomposes into two $G$-orbits: the totally real \v Silov boundary $\partial_{e}D\simeq\RR\PP^3$ and the smooth part $M$ of $\partial D$ that has Property (P). The manifold $M$ is locally CR-equivalent to the tube submanifold over the future the light cone (see (\ref{BA}) for $p=2$, $q=1$). Here $Z$ is simply-connected while the homogeneous manifold $M$ has fundamental group $\ZZ_{2}$.
\end{example}

$$
\begin{array}{l}
\hbox{Department of Mathematics}\\
\hbox{The Australian National University}\\
\hbox{Canberra, ACT 0200}\\
\hbox{Australia}\\
\hbox{E-mail: alexander.isaev@anu.edu.au}\\
\hbox{\hfill}
\end{array}\quad
\begin{array}{l}
\hbox{Mathematisches Institut}\\
\hbox{Universit\"at T\"ubingen}\\
\hbox{Auf der Morgenstelle 10}\\
\hbox{72076 T\"ubingen}\\
\hbox{Germany}\\
\hbox{E-mail: kaup@uni-tuebingen.de}
\end{array}
$$

\end{document}